\documentclass{amsart}
\usepackage{amssymb,latexsym}
\usepackage{amsfonts}

\newtheorem{thm}{Theorem}

\newtheorem{lem}{Lemma}
\newtheorem{cor}{Corollary}
\newtheorem{rem}{Remark}

\newcommand{\ba}{\begin{array}}
\newcommand{\ea}{\end{array}}

\def \qed{\cqfd}

\def \a{\alpha}

\def \b{\beta}
\def \r {\rho}
\def \p{\partial}

\def \O{\Omega}
\def \o{\omega}
\def \d{\delta}
\def \e{\epsilon}

\def \g{\gamma}
\def \G{\Gamma}
\def\n{\nabla}
\def\tn{\tilde{\nabla}}
\def \s{\sigma}

\def\th{\theta}

\def \w{\wedge}

\def \bi{\bar {i}}
\def \bj{\bar {j}}

\def \ba{\bar{\alpha}}

\def \tg{\tilde{g}}
\def \gp{g_{\phi}}
\def\into{\int_{\O}}
\def\on{\frac{\o^n}{n!}}
\def\tri{\triangle}
\def\td{\tilde}

\def\qed{\vbox{\hrule
\hbox{\vrule\hbox to 5pt{\vbox to 8pt{\vfil}\hfil}\vrule}\hrule}}

\newcommand{\beg}{\begin{eqnarray*}}
\newcommand{\begn}{\begin{eqnarray}}
\newcommand{\en}{\end{eqnarray*}}
\newcommand{\enn}{\end{eqnarray}}
\newcommand{\begu}{\begin{equation}}
\newcommand{\enu}{\end{equation}}

\title{Regularity estimates of solutions to complex Monge-Amp\`ere equations on Hermitian manifolds}

\author{Xi Zhang}
\address{Xi Zhang\\
Department of Mathematics\\
Zhejiang University, P. R. China } \email{xizhang@zju.edu.cn}
\author{Xiangwen Zhang}
\address{Xiangwen Zhang\\
Department of Mathematics and Statistics\\
McGill University, Canada} \email{xzhang@math.mcgill.ca}

\begin{document}
\maketitle
\begin{abstract}
In this paper, we obtain the Bedford-Taylor interior $C^{2}$ estimate and local Calabi $C^{3}$ estimate for the solutions to complex Monge-Amp\`ere
equations on Hermitian manifolds.
\end{abstract}

\section{Introduction}
\setcounter{equation}{0}
\par
The complex Monge-Amp\`ere equation is one of the most important partial differential equations in complex geometry. The proof of the Calabi conjecture given by S.T. Yau \cite{Y} in 1976 yields significant applications of the Monge-Amp\`ere equation in K\"ahler geometry. After that, many important geometric results, especially in K\"ahler geometry, were obtained by studying this equation. It is natural and also interesting to study the complex Monge-Amp\`ere equations in a more general form and in different geometric settings.
 \par
There are many modifications and generalizations in the existing literature. In \cite{TWY}, Tosatti-Weinkove-Yau gave a partial affirmative answer to a conjecture of Donaldson in symplectic geometry by solving (under additional curvature assumption) the complex Monge-Amp\`ere equation in an almost K\"ahler geometric setting. By studying a more general form of the Monge-Amp\`ere equation on non-K\"ahler manifolds, Fu and Yau \cite{FY} gave a solution to the Strongminger system which is motived by superstring theory. Another direction worth studying is the corresponding equation on Hermitian manifolds. In such a case the equation is not so geometric, since Hermitian metrics do not represent positive cohomology classes. On the other hand the estimates for Hermitian manifolds are more complicated than the K\"ahler case because of the non-vanishing torsion.
\par
In the eighties and nineties, some results regarding the Monge-Amp\`ere equation in the Hermitian setting were obtained by Cherrier \cite{Ch1}, \cite{Ch2} and Hanani \cite{Ha}. For next few years there was no activity on the subject until very recently, when the results were rediscovered and generalized by Guan-Li \cite{GL}. Under additional conditions they generalized the a priori estimates due to Yau \cite{Y} from the K\"ahler case and got some existence results for the solution of the complex Monge-Amp\`ere equation. At the same time, Zhang \cite{Z} independently proved similar a priori estimates in the Hermitian setting and he also considered a general form of the complex Hessian equation. Later, Tosatti-Weinkove \cite{TW1}, \cite{TW2} gave a more delicate a priori $C^2$-estimate and removed the conditions in \cite{GL}. Moreover, Dinew-Kolodziej \cite{DK} also studied the equation in the weak sense and obtained the $L^{\infty}$ estimates via suitably constructed pluripotential theory. In this paper, we want to study some other regularity properties of the complex Monge-Amp\`ere equation on Hermitian manifolds:  the Bedford-Taylor interior $C^2$-estimate and Calabi $C^3$-estimate.

\par
 The interior estimate for second order derivatives is an important and difficult topic in the study of complex Monge-Amp\`ere equation.
 It
 has many fundamental applications in complex geometric problems. In the cornerstone work of Bedford and Taylor \cite{BT}, by using the
 transitivity of the automorphism group of the unit ball $\mathbb B\subset \mathbb C^n$, they obtained the interior $C^2$-estimate for
 the following Dirichlet
 problem:
\begin{equation*}
 \begin{cases}
\det(u_{i\bj})= f \text{ in } \mathbb B,\\
u=\phi\ {\rm on}\ \p\mathbb B,
 \end{cases}
\end{equation*}
where $\phi\in \mathcal C^{1,1}(\p \mathbb B)$ and $0\leq f^{\frac1n}\in \mathcal C^{1,1}(\mathbb B)$.
\par
 Unfortunately for generic domains $\O\subset\mathbb C^n$, due to the non-transitivity of the automorphism group of $\O$, Bedford and
 Taylor's method is not applicable and the analogous estimate is still open.
Here, we exploit the method of Bedford-Taylor to study the interior estimate for the Dirichlet problem of the complex Monge-Amp\`ere equation in the unit ball in the Hermitian setting (notice that for local arguments the shape of the domain is immaterial and hence it suffices to consider balls). We consider the following equation
\begin{equation}\label{beqn}
 \begin{cases}
(\o+\sqrt{-1}\p\bar\p u)^n= f \o^n \text{ in } \mathbb B,\\
u=\phi\ {\rm on}\ \p\mathbb B,
 \end{cases}
\end{equation}
where $0\leq f^{\frac{1}{n}} \in \mathcal C^{1,1}(\mathbb B)$ and $\o$ is a smooth positive $(1, 1)$-form (not necessarily closed ) defined on $\bar{\mathbb B}$. We denote $PSH(\omega , \Omega )$ be the set of all integrable, upper semicontinuous functions satisfying $(\o+\sqrt{-1}\p\bar\p u)\geq 0$ in the current sense on the domain $\Omega$. Since $\omega $ is not necessarily K\"ahler, there are no local potentials for $\omega $, and thus Bedford-Taylor's method can not be applied directly in our case.
\begin{thm}\label{bthm2}
 Let $\mathbb B$ be the unit ball on $\mathbb C^{n}$ and $\omega $ be a smooth positive $(1, 1)$-form (not necessary closed) on
 $\bar{\mathbb B}$.   Let $u\in \mathcal C(\bar{\mathbb B}) \cap PSH(\omega , \mathbb B )$ solve the Dirichlet problem (\ref{beqn}) with
 $\phi\in \mathcal C^{1,1}(\p \mathbb B)$. Then $u\in \mathcal C^{1,1}(\mathbb B)$, and  for  arbitrary compact subset
  $B'\subset\subset\mathbb B$, there exists a constant $C$\ dependent only on $\omega $ and $\text{dist}\{ B',\p \mathbb B\}$ such that
$$||u||_{\mathcal C^{1,1}(B')}\leq C||\phi||_{\mathcal C^{1,1}(\p\mathbb B)}+ C ||f^{\frac{1}{n}}||_{\mathcal C^{1,1}(\mathbb B)}.$$
\end{thm}

\begin{rem}
Observe  that this estimate is scale and translation invariant i.e. the same constant will work if we consider the Dirichlet problem in any ball with arbitrary radius (and suitably rescaled set $B'$). \end{rem}

\medskip
As we have already mentioned, another goal of this paper is to get a local version of the $C^3$-estimate of the complex Monge-Amp\`ere equation on Hermitian manifolds. The Calabi's $C^3$-estimate for the real Monge-Amp\`ere equation was first proved by Calabi himself in \cite{Ca}. After that many mathematicians paid a lot of attention to this estimate. In Yau's celebrated work \cite{Y} about the Calabi conjecture,  he gave a detailed proof of the $C^3$-estimate for the complex Monge-Amp\`ere equation on K\"ahler manifolds, which was generalized to the Hermitian case by Cherrier \cite{Ch1}.
\par
All these $C^3$-estimates are global. However, in some situations, a local $C^3$-estimate is needed. For example  Riebesehl and Schulz \cite{SR} gave a local version of Calabi's estimate in order to study the Liouville property of Monge-Amp\`ere equations on $\mathbb{C}^n$. In a recent work by Dinew and the authors \cite{DZZ}, aimed to study the $\mathcal C^{2,\a}$ regularity of solutions to complex Monge-Amp\`ere equation, the local result in \cite{SR} also played an important role to get the optimal value of $\a$. Thus, it is also natural to generalize this local estimate to Hermitian manifolds and find some interesting applications.

\par
Let $(M, g)$ be a Hermitian manifold. We consider the following complex Monge-Amp\`ere equation
\begin{equation}\label{eqn}
(\o+\sqrt{-1}\p\bar\p \phi)^n=e^f \o^n,
\end{equation}
where $0\leq f(z)\in C^{\infty}(M)$ and $\o$ is the Hermitian form associated with the metric g.

\begin{thm}
Let $\phi(z)\in PSH(\omega , M)\cap \mathcal C^{4}(M)$ be a  solution of the Monge-Amp\`ere equation (\ref{eqn}), satisfying
\begin{eqnarray}\label{C1}
|d\phi |_{\omega }+|\partial \bar{\partial }\phi |_{\omega }\leq K.
\end{eqnarray}
Let $\O'\subset\subset \O\subset M$. Then the third derivatives of $\phi(z)$ of mixed type can be estimated in the form
\[
|\nabla_{\omega }\partial \bar{\partial }\phi |_{\omega }\leq C \quad \text{ for } z\in \O',
\]
where $C$ is a constant depending on $ K,  |d\o|_{\o},  |R|_{\o}, |\n R|_{\o}, |T|_{\o}, |\n T|_{\o},  dist (\O', \partial \O )$ and $|\n ^s f|_{\o}$ , $s=0, 1, 2, 3$. Here $\n$ is the Chern connection with respect to the Hermitian metric $\omega $, $T$ and $R$ are the torsion tensor and curvature form of $\n$.
\end{thm}

\par
 From the detailed proof in Yau's paper \cite{Y} (see also \cite{PSS}), in the K\"ahler case, we know that the quantity considered by
 Calabi
\[
S=\tg^{j\bar r}\tg^{s\bar k}\tg^{m\bar l}\phi_{j\bar k m}\phi_{\bar r s\bar l}
\]
 satisfies the following elliptic inequality:
 \begin{equation}\label{ellip1}
\tilde{\tri}S\geq -C_1 S-C_2.
\end{equation}
Here $\phi $ is a smooth solution of the equation (\ref{eqn}), $ \tilde{g}$ denotes the Hermitian metric with respect to the form $\omega_{\phi }=\omega +\sqrt{-1}\partial \bar{\partial } \phi $, $\phi_{i\bar{j}k}$ denotes the covariant derivative with respect to the Chern connection $\n$. Riebesehl and Schulz \cite{SR} used the above elliptic inequality to get the $L^p$ estimate for $S$. Then, a standard theorem for linear elliptic equations gave the $L^{\infty}$ estimate. For the Hermitian case, due to the non-vanishing torsion term, the estimates are more complicated. In \cite{Ch1}, Cherrier proved the elliptic inequality corresponding to (\ref{ellip1}) on Hermitian manifolds:
\begin{equation}\label{elliptic}
\tilde{\tri}S\geq -C_1 S^{\frac{3}{2}}-C_2,
\end{equation}
where $\tilde{\triangle }$ is the canonical Laplacian with respect to the Hermitian metric $\tilde{g}$ (i.e. $\tilde{\triangle }f=2\tg^{i\bar j}f_{i\bar j}$), positive constants $C_{1}$ and $C_{2}$ depend on $ K,   |R|_{\o}, |\n R|_{\o}, |T|_{\o}, |\n T|_{\o}$, and $|\n ^s f|_{\o}$ , $s=0, 1, 2, 3$.
\par
 By a similar method to that in \cite{SR}, we obtain the $L^p$ estimate for $S$, and then use Moser iteration to get the $L^{\infty}$
 estimate.

\par
The estimates obtained in this paper should be useful for the study of problems on Hermitian manifolds. As a simple application, following the lines of \cite{DZZ}, one has the following corollary:
\begin{cor} Let $\Omega$ be a domain in $\mathbb C^n$ and $\omega$ be a Hermitian form defined on $\Omega$. Let $\phi(z)\in PSH(\omega , \Omega)\cap \mathcal C^{1,1}(\Omega)$ be a  solution of the Monge-Amp\`ere equation
$$ (\o+\sqrt{-1}\p\bar\p \phi)^n=e^f \o^n.$$
 Suppose that $f$ is strictly positive and $f\in\mathcal C^{\alpha}(\Omega)$ for some $0<\alpha<1$. Then $\phi\in\mathcal
 C^{2,\alpha}(\Omega)$.
\end{cor}

\medskip
\par
The paper is organized as follows. In section 2, we prove the interior $C^2$-estimate for the complex Monge-Amp\`ere equation. The proof for Calabi's $C^3$-estimate is given in section 3. In the appendix, we give a new proof of (\ref{elliptic}) which follows the idea in Phong-Sesum-Sturm \cite{PSS}, where the authors gave a simpler proof of Calabi's estimate on K\"ahler manifolds.

\vspace{1cm}
\section{Proof of the interior estimates}
\par

\medskip
\par
In the proof of interior $C^2$-estimates, the comparison theorem will play the key role. Following the same idea as in \cite{CKNS}, it's easy to see that the comparison theorem is still true for the complex Monge-Amp\`ere equation on Hermitian manifold $(M, \omega)$.

\begin{lem}
  Let $\O\subset$ M be a bounded set and  $u, v \in C^{\infty}(\bar \O)$, with $\o+\sqrt{-1}\p\bar\p u\geq 0$, $\o+\sqrt{-1}\p\bar\p v> 0$  be such that
\begin{equation*}
(\o+\sqrt{-1}\p\bar\p v)^n\geq (\o+\sqrt{-1}\p\bar\p u)^n
\end{equation*}
and
\begin{equation*}
v\leq u  \quad \text{ on } \p \O,
\end{equation*}
then $v\leq u$ in $\bar \O.$
\end{lem}

\medskip

\noindent{\bf Proof of Theorem 1:}
\par
As mentioned above, we will follow the idea of Bedford and Taylor from \cite{BT}.
For $a\in B^n$, let $T_a\in \text{Aut} (B^n)$ be defined by
\begin{equation*}
T_a(z)=\G(a)\frac{z-a}{1-\bar a^t z},
\end{equation*}
where $\G(a)=\frac{a^t \bar a}{1-v(a)}-v(a)I$ and $v(a)=\sqrt{1-|a|^2}.$
\par
Note that $T_a(a)=0, T_{-a}=T_a^{-1}$, and $T_a(z)$ is holomorphic in $z$, and a smooth function in $a\in B^n$. For any $a\in B(0, 1-\eta)=\{a: |a|<1-\eta\}$ set
\begin{equation*}
L(a, h, z)=T^{-1}_{a+h}T_a(z)
\end{equation*}
and
\begin{eqnarray*}
&&U(a, h, z)=L^*_1 u(z), \quad U(a, -h, z)=L^*_2 u(z), \\
&&\Phi(a, h, z)=L^*_1 \phi(z), \quad \Phi(a, -h, z)=L^*_2 \phi(z), \quad \text{ for } z\in \p B^n.
\end{eqnarray*}
where $L^*_i$ means the pull-back of $L_i \text{ for } i=1, 2$ and $L_1=L(a, h, z), L_2=L(a, -h, z)$. Since $U(a, h, z)=\Phi(a, h, z)$ for $z\in \p B^n$, it follows that $U\in C^{1,1}(B(0, 1-\eta)\times B(0, \eta)\times \p B^n)$. Consequently, for a suitable constant $K_1$, depending on $\eta>0$, we have
\begin{equation}\label{b2}
\frac{1}{2}(U(a, h, z)+U(a, -h, z))-K_1 |h|^2\leq U(a, 0, z)=\phi(z)
\end{equation}
for all $|a|\leq 1-\eta, |h|\leq \frac{1}{2}\eta,\text{ and } z\in \p B^n$. If it can be shown that $v(a, h, z)$ satisfies
\begin{equation}\label{b3}
(\o+\sqrt{-1}\p\bar\p v)^n\geq f(z) \o^n,
\end{equation}
where
\begin{equation}\label{bv}
v(a, h, z)=\frac{1}{2}\Big[U(a, h, z)+U(a, -h, z)\Big]-K_1 |h|^2+K_2(|z|^2-1)|h|^2,
\end{equation}
then it follows from the comparison theorem in the Hermitian case that $v(a, h, z)\leq u(z)$. Thus, if we set $a=z$, we conclude that
\begin{equation*}
\frac{1}{2}[u(z+h)+u(z-h)]\leq u(z)+(K_1+K_2)|h|^2
\end{equation*}
which would prove the theorem.
\par
Let now
\begin{eqnarray}\label{b4}
F(\o+\sqrt{-1}\p\bar\p v) &=&\Big(\frac{(\o+\sqrt{-1}\p\bar\p v)^n}{(\sqrt{-1})^n dz^1\wedge d\bar z^1\wedge\cdots\wedge dz^n\wedge d\bar z^n}\Big)^{\frac{1}{n}}\\\nonumber
&=&\Big(\det(g_{i\bj}+v_{i\bj})\Big)^{\frac{1}{n}},
\end{eqnarray}
where $g_{i\bj}$ is the local expression of $\o$ under the standard coordinate $\{z_{i}\}_{i=1}^{n}$ in $\mathbb C^{n}$.
\par
By the concavity of $F$, we have
\begin{eqnarray}\label{b5}
F(\o+\sqrt{-1}\p\bar\p v)&=& F\Big(\o+\frac{\sqrt{-1}}{2}(\p\bar \p L_1^* u +\p\bar \p L_2^* u+2K_2|h|^2 \p\bar \p |z|^2)\Big)\\\nonumber
&=& F\Big(\frac{1}{2}(\o-L_1^*\o)+\frac{1}{2}(\o-L_2^*\o)+K_2|h|^2\sqrt{-1}\p\bar\p|z|^2\\\nonumber
&&+\frac{1}{2}(L_1^*\o+\sqrt{-1}\p\bar\p L_1^* u)+\frac{1}{2}(L_2^*\o+\sqrt{-1}\p\bar\p L_2^* u)\Big)\\\nonumber
&\geq& \frac{1}{2}F(L_1^*\o+\sqrt{-1}\p\bar\p L_1^* u)+ \frac{1}{2}F(L_2^*\o+\sqrt{-1}\p\bar\p L_2^* u)\\\nonumber
&&+\frac{1}{2}F\Big((\o-L_1^*\o)+(\o-L_2^*\o)+2K_2|h|^2\sqrt{-1}\p\bar\p|z|^2\Big).
\end{eqnarray}
Since the Hermitian metric $\o$ is smooth, one can find $K_2$ large enough, such that
\begin{equation}\label{bw}
(\o-L_1^*\o)+(\o-L_2^*\o)+K_2|h|^2\sqrt{-1}\p\bar\p|z|^2\geq 0.
\end{equation}
On the other hand, since $L(a, h, z)$ is holomorphic in $z$, it follows from the equation (\ref{beqn}) that
\begin{eqnarray}\label{b6}
F(L_1^*\o+\sqrt{-1}\p\bar\p L_1^* u)&=&F(L_1^*(\o+\sqrt{-1}\p\bar\p u))\\\nonumber
&=&\Big(\frac{L_1^*(\o+\sqrt{-1}\p\bar\p u)^n}{(\sqrt{-1})^n dz^1\wedge d\bar z^1\wedge\cdots\wedge dz^n\wedge d\bar z^n}\Big)^{\frac{1}{n}}\\\nonumber
&=& \Big(\frac{L_1^*(f(z) \o^n)}{(\sqrt{-1})^n dz^1\wedge d\bar z^1\wedge\cdots\wedge dz^n\wedge d\bar z^n}\Big)^{\frac{1}{n}}\\\nonumber
&=&F(L_1^*(f^{\frac{1}{n}}  \o))=L_{1}^{\ast}(f^{\frac{1}{n}})F(L_1^*(  \o)).
\end{eqnarray}
Similarly, we can get
\[
F(L_2^*\o+\sqrt{-1}\p\bar\p L_2^* u)=F(L_2^*(f^{\frac{1}{n}}  \o))=L_{2}^{\ast}(f^{\frac{1}{n}})F(L_2^*(  \o)).
\]
Thus,
\begin{eqnarray}\label{b7}
\nonumber F(\o+\sqrt{-1}\p\bar\p v)&\geq& \frac{1}{2} \Big(F(L_1^*(f^{\frac{1}{n}}  \o))+ F(L_2^*(f^{\frac{1}{n}}  \o))\Big)+\frac{1}{2}F(K_2|h|^2\sqrt{-1}\p\bar\p|z|^2)\\
&=&F(f^{\frac{1}{n}}  \o)+\frac{1}{2} \Big(F(L_1^*(f^{\frac{1}{n}}  \o))+ F(L_2^*(f^{\frac{1}{n}}  \o))-2F(f^{\frac{1}{n}}  \o)\Big)\\\nonumber
&&+\frac{1}{2}F(K_2|h|^2\sqrt{-1}\p\bar\p|z|^2).\\\nonumber
\end{eqnarray}
Again, since $\o$ is smooth and $f^{1/n}\in C^{1,1}$, choosing $K_2$ large enough, we have
\begin{eqnarray}
F(L_1^*(f^{\frac{1}{n}}  \o))+ F(L_2^*(f^{\frac{1}{n}}  \o))-2F(f^{\frac{1}{n}}  \o)\leq F(K_2|h|^2\sqrt{-1}\p\bar\p|z|^2).
\end{eqnarray}
Finally, we obtain
\begin{equation}\label{b8}
F(\o+\sqrt{-1}\p\bar\p v)\geq F(f^{\frac{1}{n}}  \o),
\end{equation}
and thus, the inequality (\ref{b3}) follows.

\medskip

\vspace{1cm}
\section{Proof of the Calabi estimate}

\par
Let $(M, J, \o)$ be a Hermitian manifold and $\n$ denote the Chern connection with respect to the metric $\o$. Let locally $\o=\sqrt{-1}g_{i\bj}dz^i\w dz^{\bj}$, then the local formula for the connection 1-form reads $\th=\p g\cdot g^{-1}$. We also denote
\[
\th_{\a}=\p_{\a}g\cdot g^{-1}, \quad  \th^{\g}_{\a\b}=\frac{\p g_{\b\bar \d}}{\p z^{\a}}g^{\g\bar \d}.
\]
The torsion tensor of $\n$ is defined by
\begin{eqnarray*}
T(\frac{\p}{\p z^{\a}}, \frac{\p}{\p z^{\b}})&=&\n_{\frac{\p}{\p z^{\a}}}\frac{\p}{\p z^{\b}}-\n_{\frac{\p}{\p z^{\b}}}\frac{\p}{\p z^{\a}}-[\frac{\p}{\p z^{\a}}, \frac{\p}{\p z^{\b}}]\\
&=&\Big(\frac{\p g_{\b \bar \d}}{\p z^{\a}}-\frac{\p g_{\a\bar\d}}{\p z^{\b}}\Big)g^{\g\bar \d}.
\end{eqnarray*}
Notice that $T=0\iff \o$ is K\"ahler (and $\n$ is the Levi-Civita connection on M).
\par
The curvature form of $\n$ is defined by $R=\bar \p \th= d\th-\th\w\th=\bar \p (\p g\cdot g^{-1})$. In local coordinates, we have
\begin{eqnarray*}
&&R^j_{i\a\bar\b}=-\bar\p_{\b}(\p_{\a}g\cdot g^{-1})^j_i
=-g^{j\bar k}\frac{\p^2 g_{i\bar k}}{\p z^{\a}\p \bar z^{\b}}
+\frac{\p g_{i\bar k}}{\p z^{\a}}g^{j\bar s}
\frac{\p g_{t\bar s}}{\p \bar z^{\b}}g^{t\bar k},\\
&& R_{i\bj \a\bar \b}=g_{k\bj}R^k_{i\a\bar \b}.
\end{eqnarray*}
Note that $R^{(2, 0)}=R^{(0, 2)}=0 \text{ and } T^{(1, 1)}=0$, since the almost complex structure $J$ is integrable and $\n$ is the Chern connection.

\medskip
\noindent{\bf Proof of Theorem 2:}
\par
By the assumption (\ref{C1})  for the solution of equation (\ref{eqn}), we know that
\[
\frac{1}{\lambda} g\leq g_{\phi}\leq \lambda g \quad \text{ for some constant }\lambda ,
\]
where $\lambda $ depends only on $K$ and $\|f\|_{C^{0}}$, and $g_{\phi }$ denotes the Hermitian metric with respect to the form $\omega_{\phi }=\omega +\sqrt{-1}\partial \bar{\partial } \phi $. Thus,
\begin{eqnarray}\label{2}
S=(\gp)^{j\bar r}(\gp)^{s\bar k}(\gp)^{m\bar l}\phi_{j\bar k m}\phi_{\bar r s \bar l}
\leq \lambda (\gp)^{j\bar r}(\gp)^{s\bar k}g^{m\bar l}\phi_{j\bar k m}\phi_{\bar r s \bar l}.
\end{eqnarray}
On the other hand, we have
\begin{eqnarray*}
\gp^{j\bar k}g^{m\bar l}\phi_{j\bar k m\bar l}&=&\Big (\gp^{j\bar k}g^{m\bar l}\phi_{j\bar k m}\Big)_{\bar l}-(\gp^{j\bar k})_{\bar l} g^{m\bar l}\phi_{j\bar k m}\\\nonumber
&=& g^{m\bar l}f_{m\bar l}+\gp^{j\bar s}\phi_{t\bar s \bar l}\gp^{t\bar k}g^{m \bar l}\phi_{j\bar k m},
\end{eqnarray*}
where we used the equation (\ref{eqn}) in the last equality above. Thus
\begin{equation}\label{3}
S\leq \lambda\Big[\gp^{j\bar k}g^{m \bar l}\phi_{j\bar k m\bar l}-\triangle f \Big].
\end{equation}
Notice that $\gp^{j\bar k}g^{m \bar l}\phi_{j\bar k m\bar l}=\wedge_{\gp}(g^{m\bar l}\nabla_{\bar l}\nabla_m(\sqrt{-1}\p\bar \p \phi))$ is a globally defined quantity.
Therefore we can estimate for every sufficiently large exponents $\r, \s,$ and every nonnegative test function $\eta(z)\in C^1_0(\O)$:
\begin{equation}\label{4}
\into S^{\s}\eta^{p+1}\on\leq \lambda \into S^{\s-1}\eta^{p+1}[\gp^{j\bar k}g^{m\bar l}\phi_{j\bar k m\bar l}-\triangle f]\on.
\end{equation}
Now, using the following identity:
\begin{eqnarray*}
\phi_{j\bar k m\bar l}&=&\phi_{j\bar k\bar l m}+\phi_{s\bar k}R^s_{j m \bar l}-\phi_{j\bar t}R^{\bar t}_{\bar k m\bar l}\\
&=&\phi_{j\bar l m \bar k}+\phi_{s\bar l}R^s_{j m \bar k}+\phi_{s\bar k}R^s _{j m \bar l}-\phi_{j\bar t}R^{\bar t}_{\bar l m\bar k}-\phi_{j\bar k}R^{\bar t}_{\bar k m\bar l}\\
&=& \phi_{m\bar l j\bar k}+ C_1,
\end{eqnarray*}
where $C_1$ is a constant depending on $K$ and $|R|_{\o}$.
Therefore, we have
\begin{eqnarray}\label{5}
\into S^{\s}\eta^{p+1}\on&\leq& \lambda (\into S^{\s-1}\eta^{p+1}\gp^{j\bar k}g^{m\bar l}\phi_{m\bar l j\bar k}\on+\\\nonumber &&\into S^{\s-1}\eta^{p+1}(C_1-\triangle f)\on)\\\nonumber
&\leq& \lambda \into S^{\s-1}\eta^{p+1}\gp^{j\bar k}(\triangle \phi)_{j\bar k}\on +C_2 \into S^{\s-1}\eta^{p+1}\on,
\end{eqnarray}
where $C_2$ is a constant depending on $C_1 \text{ and } \triangle f$.

\par
Now, using integration by parts, it is easy to see that
\begin{eqnarray*}
&&\into S^{\s-1}\eta^{p+1}\gp^{j\bar k}(\triangle \phi)_{j\bar k}\on\\
&=& \into e^{-f}S^{\s-1}\eta^{p+1}\gp^{j\bar k}(\triangle \phi)_{j\bar k}\frac{\o_{\phi}^{n}}{n!}\\
&=&\into e^{-f}S^{\s-1}\eta^{p+1}\sqrt{-1}\p\bar\p(\triangle \phi)\wedge\frac{\o_{\phi}^{n-1}}{(n-1)!}\\
&=&\into \sqrt{-1}d(e^{-f}S^{\s-1}\eta^{p+1}\bar\p(\triangle\phi))\wedge \frac{\o_{\phi}^{n-1}}{(n-1)!}\\&&-\into \sqrt{-1}d(e^{-f}S^{\s-1}\eta^{p+1})\wedge\bar\p(\triangle\phi)\wedge \frac{\o_{\phi}^{n-1}}{(n-1)!}\\
&=:&I-II.
\end{eqnarray*}
Next, we will estimate $I$ and $II$. First,
\begin{eqnarray}\label{6}
I&=&\into \sqrt{-1}d(e^{-f}S^{\s-1}\eta^{p+1}\bar\p(\triangle\phi))\wedge \frac{\o_{\phi}^{n-1}}{(n-1)!}\\\nonumber
&=&-\into\sqrt{-1}e^{-f}S^{\s-1}\eta^{p+1}\bar\p(\triangle \phi)\wedge d \o_{\phi}\wedge \frac{\o_{\phi}^{n-2}}{(n-2)!}.
\end{eqnarray}
By the equivalence of two forms $\o \text{ and } \o_{\phi}$ (i.e., the assumption (\ref{C1}) on $\phi$), we know
\begin{eqnarray}\label{7}
\bar\p(\triangle \phi)\wedge d\o_{\phi}\wedge\frac{\o_{\phi}^{n-2}}{(n-2)!}&=& \bar\p(\triangle \phi)\wedge d\o\wedge \frac{\o_{\phi}^{n-2}}{(n-2)!}\\\nonumber
&\leq&C_3 |\bar\p (\triangle \phi)|_{g_{\phi}}|d\o|_{g_{\phi}} \on\\\nonumber
&\leq& C_4 S^{\frac{1}{2}}\on,
\end{eqnarray}
where $C_4$ is a constant depending on $|d\o|_g, \|f\|_{C^{0}} \text{ and } K$ (for the justification of the last inequality we refer to the formula of $S$ given in the appendix). This estimate yields
\begin{equation}\label{8}
I\leq C_5 \into S^{\s-\frac{1}{2}}\eta^{p+1}\on
\end{equation}
for some constant $C_5$ dependent on $\o, \|f\|_{C^{0}} \text{ and } K$.

Let us now estimate the second term:
\begin{eqnarray}\label{9}
II&=& \into \sqrt{-1}d(e^{-f})S^{\s-1}\eta^{p+1}\w \bar\p(\tri \phi)\w\frac{\o_{\phi}^{n-1}}{(n-1)!}\\\nonumber
&&+(\s-1)\into \sqrt{-1}e^{-f}S^{\s-2}\eta^{p+1}dS\w\bar\p(\tri \phi)\w\frac{\o_{\phi}^{n-1}}{(n-1)!}\\\nonumber
&&+(p+1)\into \sqrt{-1}e^{-f}S^{\s-1}\eta^p d\eta\w\bar\p(\tri\phi)\w\frac{\o_{\phi}^{n-1}}{(n-1)!}\\\nonumber
&\leq&C_6\Big(\into S^{\s-\frac{1}{2}}\eta^{p+1}\on +(\s-1)\into S^{\s-\frac{3}{2}}|\nabla S|\eta^{p+1}\on\\\nonumber
&&+(p+1)\into S^{\s-\frac{1}{2}}\eta^p|\nabla \eta|\on\Big),
\end{eqnarray}
where $C_6$ is a constant depending on $\|f\|_{C^{1}(\o )} \text{ and } K$.
\par
By the estimates (\ref{8}), (\ref{9}) and using Cauchy's inequality
\begin{equation*}
(\s-1)\eta^{p+1}S^{\s-\frac{3}{2}}|\nabla S|\leq \frac{(\s-1)^2}{4\e}\eta^{p+1}S^{\s-3}|\nabla S|^2+\e \eta^{p+1}S^{\s}
\end{equation*}
we have, for $\e>0$ small enough,
\begin{eqnarray}\label{10}
\into S^{\s}\eta^{p+1}\on&\leq& C_7\Big((\s-1)^2\into S^{\s-3}|\nabla S|^2\eta^{p+1}\on +\into S^{\s-1}\eta^{p+1}\on\\\nonumber
&&+(p+1)\into S^{\s-\frac{1}{2}}\eta^p|\nabla\eta|\on +\into S^{\s-\frac{1}{2}}\eta^{p+1}\on\Big),
\end{eqnarray}
where $C_7$ is a constant depending on $|d\o |_{\o}, |R|_{\o}, K, \|f\|_{C^{1}(\o )} \text{ and }\tri f$.
\medskip
\par
Now we are in the place to use the elliptic inequality (\ref{elliptic}) in the introduction. Recall that
\begin{equation}\label{11}
\tri_{\phi}S\geq -CS^{\frac{3}{2}}-C_0.
\end{equation}
Multiplying by $S^{\s-2}\eta^{p+1}$ on both sides of the above inequality and integrating over $\O$, we have
\begin{equation}\label{12}
-C\into S^{\s-\frac{1}{2}}\eta^{p+1}\on -C_0\into S^{\s-2}\eta^{p+1}\on\leq \into S^{\s-2}\eta^{p+1}\tri_{\phi}S\on.
\end{equation}
The right hand side of above inequality can be estimated as follows
\begin{eqnarray*}
&&\into S^{\s-2}\eta^{p+1}\tri_{\phi}S\on\\
&=&\into e^{-f}S^{\s-2}\eta^{p+1}\sqrt{-1}\p\bar\p S\w \frac{\o_{\phi}^{n-1}}{(n-1)!}\\
&=&\into\sqrt{-1}d(e^{-f}S^{\s-2}\eta^{p+1}\bar\p S)\w\frac{\o_{\phi}^{n-1}}{(n-1)!}-\into \sqrt{-1}d(e^{-f}S^{\s-2}\eta^{p+1})\w\bar\p S\w \frac{\o_{\phi}^{n-1}}{(n-1)!}\\
&=&-\into \sqrt{-1}e^{-f}S^{\s-2}\eta^{p+1}\bar\p S\w d\o\w\frac{\o_{\phi}^{n-2}}{(n-2)!}-\sqrt{-1}\into d(e^{-f})S^{\s-2}\eta^{p+1}\w\bar\p S\w\frac{\o_{\phi}^{n-1}}{(n-1)!}\\
&&-(\s-2)\into \sqrt{-1}e^{-f}S^{\s-3}\eta^{p+1}\p S\w\bar\p S\w\frac{\o_{\phi}^{n-1}}{(n-1)!}\\
&&-(p+1)\into\sqrt{-1}e^{-f}S^{\s-2}\eta^p\p \eta\w\bar\p S\w\frac{\o_{\phi}^{n-1}}{(n-1)!}\\
&\leq&-C_8(\s-2)\into S^{\s-3}\eta^{p+1}|\nabla S|^2\on+C_9\into S^{\s-2}\eta^{p+1}|\nabla S|\on \\
&&+C_9(p+1)\into S^{\s-2}\eta^p |\nabla \eta||\nabla S|\on.
\end{eqnarray*}
From this, we obtain,
\begin{eqnarray}\label{13}
&&(\s-2)\into S^{\s-3}\eta^{p+1}|\nabla S|^2\on\\\nonumber
&\leq& C_{10}\Big((p+1)\into S^{\s-2}\eta^p |\nabla \eta||\nabla S|\on+\into S^{\s-2}\eta^{p+1}|\nabla S|\on \\\nonumber
&&+\into S^{\s-\frac{1}{2}}\eta^{p+1}\on +\into S^{\s-2}\eta^{p+1}\on\Big).
\end{eqnarray}
Now, by Cauchy's inequality again,
\begin{eqnarray*}
S^{\s-2}\eta^{p+1}|\nabla S|&\leq& \e |\nabla S|^2 S^{\s-3}\eta^{p+1}+\frac{1}{4\e}\eta^{p+1}S^{\s-1}\\
(p+1)S^{\s-2}\eta^p|\nabla \eta||\nabla S|&\leq &\e |\nabla S|^2 S^{\s-3}\eta^{p+1}+\frac{(p+1)^2}{4\e}\eta^{p-1}S^{\s-1}|\nabla \eta|^2.
\end{eqnarray*}
 These two inequalities, together with (\ref{13}) and (\ref{10}) yield
\begin{eqnarray}\label{14}
&&\into S^{\s}\eta^{p+1}\on\\\nonumber
 &\leq&C_{11} \s^2(p+1)^2\Big(\into S^{\s-\frac{1}{2}}\eta^{p+1}\on +\into S^{\s-\frac{1}{2}}\eta^{p+1}\on \\\nonumber
&&+\into S^{\s-\frac{1}{2}}\eta^p|\nabla\eta|\on +\into S^{\s-2}\eta^{p+1}\on \into S^{\s-1}\eta^{p-1}|\nabla \eta|^2\on\Big)
\end{eqnarray}
for $p\geq 2, \s\geq 4$.
\medskip
\par
Now, let $B_{R_0}(z)\subset\subset\O$ be a ball, and let $0<R\leq r<t\leq R_0, R_0-R\leq 1$. By choosing an appropriate testing function $\eta(z)$, with $0\leq\eta\leq 1, \eta|_{B_r}=1, \eta|_{M/B_t}=0, |\nabla \eta|\leq \frac{C}{t-r}$, and puting $p=\s-1$, we conclude that
\begin{eqnarray}\label{15}
\nonumber\int_{B_t(z)}(S\eta)^{\s}\on
&\leq&C_{12}\s^4\int_{B_t(z)}\Big\{\frac{1}{(t-r)^2}(S\eta)^{\s-2}S\\&&+\frac{1}{t-r}(S\eta)^{\s-1}S^{\frac{1}{2}}+(S\eta)^{\s-\frac{1}{2}}\eta^{\frac{1}{2}}+(S\eta)^{\s-1}\eta+(S\eta)^{\s-2}\eta^2\Big\}\on
\end{eqnarray}
By Young's inequality
\begin{equation*}
ab\leq \e\frac{a^{\a}}{\a}+\frac{1}{\e^{\b/\a}}\frac{b^{\b}}{\b},\quad \text{ for } \e>0, \frac{1}{\a}+\frac{1}{\b}=1.
\end{equation*}
It follows that,
\begin{eqnarray*}
\frac{1}{t-r}(S\eta)^{\s-1}S^{\frac{1}{2}}&\leq&\frac{\e}{\frac{\s}{\s-1}}\Big((S\eta)^{\s-1}\Big)^{\frac{\s}{\s-1}}+\frac{1}{\e^{\s-1}\s}\Big(\frac{1}{t-r}S^{\frac{1}{2}}\Big)^{\s};  \a=\frac{\s}{\s-1}, \b=\s\\
\frac{1}{(t-r)^2}(S\eta)^{\s-2}S&\leq&\frac{\e}{\frac{\s}{\s-2}}\Big((S\eta)^{\s-2}\Big)^{\frac{\s}{\s-2}}+\frac{1}{\e^{\frac{\s-2}{2}}\frac{\s}{2}}\Big(\frac{1}{(t-r)^2}S\Big)^{\frac{\s}{2}};  \a=\frac{\s}{\s-2}, \b=\frac{\s}{2}\\
(S\eta)^{\s-2}&\leq&\frac{\e}{\frac{\s}{\s-4}}\Big((S\eta)^{\s-4}\Big)^{\frac{\s}{\s-4}}+\frac{1}{\e^{\frac{\s-4}{4}}\frac{\s}{4}}\Big( (S\eta)^2\Big)^{\frac{\s}{4}};  \a=\frac{\s}{\s-4}, \b=\frac{\s}{4}\\
(S\eta)^{\s-1}&\leq&\frac{\e}{\frac{\s}{\s-2}}\Big((S\eta)^{\s-2}\Big)^{\frac{\s}{\s-2}}+\frac{1}{\e^{\frac{\s-2}{2}}\frac{\s}{2}}\Big( S\eta\Big)^{\frac{\s}{2}};  \a=\frac{\s}{\s-2}, \b=\frac{\s}{2}\\
(S\eta)^{\s-\frac{1}{2}}&\leq&\frac{\e}{\frac{\s}{\s-1}}\Big((S\eta)^{\s-1}\Big)^{\frac{\s}{\s-1}}+\frac{1}{\e^{\s-1}\s}\Big( (S\eta)^{\frac{1}{2}}\Big)^{\s};  \a=\frac{\s}{\s-1}, \b=\s.
\end{eqnarray*}
All the above inequalities combined with (\ref{15}), lead to
\begin{eqnarray}\label{16}
\int_{B_r(z)}S^{\s}\on&\leq& C_{13}B(\e)^{\s}\Big(\frac{1}{(t-r)^{\s}}+\frac{1}{(t-r)^{\frac{\s}{2}}}+1\Big)\int_{B_t(z)}S^{\frac{\s}{2}}\on\\\nonumber
&\leq&C_{13}\frac{B(\e)^{\s}t^n}{(t-r)^{\s}}\Big(\int_{B_t(z)}S^{\s}\on\Big)^{\frac{1}{2}},
\end{eqnarray}
where $B(\e)$ is a constant depending on $\e$ which comes from the coefficients in the Young's inequalities above.

\medskip
\par
Now we can apply the Meyers' lemma:
\begin{lem}[\cite{M}]
If $u=u(x)$ is a nonnegative, non-decreasing continuous function in the interval $[0, d)$, which satisfies the functional inequality:
\begin{equation*}
u(s)\leq\frac{c}{r-s}\Big(u(r)\Big)^{1-\a}, \quad \text{ for any } 0\leq s<r<d,
\end{equation*}
with
$\a \text{ and } c$ being constants $(0<\a<1)$, then
\begin{equation*}
u(0)\leq \Big(\frac{2^{\a+1}c}{(2^{\a}-1)d}\Big)^{\frac{1}{\a}}.
\end{equation*}
\end{lem}

\medskip
\par
Using (\ref{16}) and applying the Meyers' lemma with $ d=R_0-R, s=r-R \text{ and } \phi(s)=\Big(\int_{B_{R+s}(z)}S^{\s}\on\Big)^{\frac{1}{\s}}$, one can obtain
\begin{equation*}
\phi(0)\leq \frac{C^{\frac{1}{\s}}B(\e)R_0^{\frac{1}{\s}}}{(R_0-R)^2},
\end{equation*}
and thus
\begin{equation}\label{17}
\Big(\int_{B_R(z)}S^{\s}\on\Big)^{\frac{1}{\s}}\leq \frac{(CR_0)^{\frac{1}{\s}}}{(R_0-R)^2}B(\e).
\end{equation}
From this, we obtain the $L^p$ estimate of $S$ for arbitrary $p$. However, by tracking the constant $B(\e)$, one can find that $B(\e)\sim \s^4$. Thus, we cannot get the estimate for $\sup_{\O} S$ by letting $\s\longrightarrow \infty$. We should instead use the standard Moser iteration to finish the $L^{\infty}$ estimate for $S$.

\medskip
\par
Recall that by inequality (\ref{13}) we have
\begin{eqnarray*}
&&(\s-2)\into S^{\s-3}\eta^{p+1}|\nabla S|^2\on\\
&\leq& C_{10}\Big((p+1)\into S^{\s-2}\eta^p |\nabla \eta||\nabla S|\on+\into S^{\s-2}\eta^{p+1}|\nabla S|\on\\
&&
+\into S^{\s-\frac{1}{2}}\eta^{p+1}\on +\into S^{\s-2}\eta^{p+1}\on\Big).
\end{eqnarray*}
Coupling this with Young inequalities
\begin{eqnarray*}
S^{\s-2}\eta^{p+1}|\nabla S| &\leq& \e |\nabla S|^2S^{\s-3}\eta^{p+1}+\frac{1}{4\e}\eta^{p+1}S^{\s-1},\\
(p+1)S^{\s-2}\eta^p|\nabla \eta||\nabla S|&\leq& \e |\nabla S|^2 S^{\s-3}\eta^{p+1}+\frac{(p+1)^2}{4\e}\eta^{p-1}S^{\s-1}|\nabla S|^2
\end{eqnarray*}
we have
\begin{eqnarray}\label{18}
&&(\s-2)\into S^{\s-3}\eta^{p+1}|\nabla S|^2\on\\\nonumber
 &\leq& C_{14}\into \frac{(p+1)^2}{\s-2}\eta^{p-1}S^{\s-1}|\nabla
 \eta|^2+\frac{1}{\s-2}S^{\s-1}\eta^{p+1}+S^{\s-\frac{1}{2}}\eta^{p+1}+S^{\s-2}\eta^{p+1}\on.
\end{eqnarray}
Let now $q=\s-1\geq 2, \text{ and } p=1$, then one obtains
\begin{eqnarray}\label{19}
&&\into S^{q-2}\eta^{2}|\nabla S|^2\on\\\nonumber
 &\leq& C_{15}\into \frac{1}{(q-1)^2}S^{q}|\nabla
 \eta|^2+\frac{1}{(q-1)^2}S^{q}\eta^{2}+\frac{1}{q-1}S^{q+\frac{1}{2}}\eta^{2}+\frac{1}{q-1}S^{q-1}\eta^{2}\on.
\end{eqnarray}
By the Sobolev inequality
\begin{equation*}
\Big(\into v^{\frac{2m}{m-1}}\on\Big)^{\frac{m-1}{2m}}\leq C\Big(\into |\nabla v|^2\on\Big)^{\frac{1}{2}}+C\Big(\into v^2\on\Big)^{\frac{1}{2}}
\end{equation*}
applied to  $v=\eta S^{\frac{q}{2}}$, we conclude that
\begin{eqnarray}\label{20}
&&\Big(\into (\eta S^{\frac{q}{2}})^{\frac{2m}{m-1}}\on\Big)^{\frac{m-1}{2m}}\\\nonumber
&\leq& C_{16} \Big[\Big(\into |\nabla (\eta S^{\frac{q}{2}})|^2\on\Big)^{\frac{1}{2}}+\Big(\into (\eta S^{\frac{q}{2}})^2\on\Big)^{\frac{1}{2}}\Big]\\\nonumber
&\leq& C_{17}\Big[\Big(\into S^q |\nabla \eta|^2 +(\frac{q}{2})^2S^{q-2}\eta^2 |\nabla S|^2\on\Big)^{\frac{1}{2}}+\Big(\into \eta^2 S^q\on\Big)^{\frac{1}{2}}\Big].
\end{eqnarray}
Using the inequality (\ref{19}), we have
\begin{eqnarray}\label{21}
&&\Big(\into (\eta^2 S^q)^{\frac{m}{m-1}}\on\Big)^{\frac{m-1}{m}}\\\nonumber
&\leq& C_{18}\into\Big( |\nabla \eta|^2 S^q +\eta^2S^q +\frac{q^2}{(q-1)^2}S^q |\nabla \eta|^2 +\frac{q^2}{(q-1)^2}S^q \eta^2\\\nonumber
&&+\frac{q^2}{q-1}S^{q+\frac{1}{2}}\eta^2+\frac{q^2}{q-1}S^{q-1}\eta^2 \Big)\on
\end{eqnarray}
for any $q>4$.
\par
Again, let $B_{R_0}(z)\subset\subset\O$ be a ball, and let $0<R\leq r_1<r_2\leq R_0, R_0-R\leq 1$. By choosing an appropriate testing function $\eta(z)$, with $0\leq\eta\leq 1, \eta|_{B_{r_1}}=1, \eta|_{M/B_{r_2}}=0, |\nabla \eta|\leq \frac{C}{r_2-r_1}$, we conclude that
\begin{eqnarray}\label{22}
&&\Big(\int_{B_{r_1}(z)} S^{q\frac{m}{m-1}}\on\Big)^{\frac{m-1}{m}}\\\nonumber
&\leq& C_{19} \int_{B_{r_2}(z)}\Big((1+\frac{q^2}{(q-1)^2})(\frac{1}{(r_2-r_1)^2}+1)S^q +\frac{q^2}{q-1}S^{q+\frac{1}{2}}+\frac{q^2}{q-1}S^{q-1}\Big)\on\\\nonumber
&\leq & qC_{20}(\frac{1}{(r_2-r_1)^2}+1)\int_{B_{r_2}(z)}(S^q+S^{q-1}+S^{q+\frac{1}{2}})\on\\\nonumber
&\leq &qC_{21}(\frac{1}{(r_2-r_1)^2}+1)\int_{B_{r_2}(z)}S^{q+\frac{1}{2}}\on.
\end{eqnarray}
Thus,
\begin{equation}\label{23}
||S||_{L^{\frac{qm}{m-1}}(B_{r_1}(z))}\leq \Big[C q(\frac{1}{(r_2-r_1)^2}+1)\Big]^{\frac{1}{q}}||S||^{\frac{q+\frac{1}{2}}{q}}_{L^{q+\frac{1}{2}}(B_{r_2}(z))}
\end{equation}
 for any $0<R\leq r_1<r_2\leq R_0$.
\par
Let $\frac{q_k m}{m-1}=q_{k+1}+\frac{1}{2}$ and $r_k=R+(R_0-R)2^{-k}$. Then,
\[
q_k=\Big(\frac{m}{m-1}\Big)^k+\frac{m-1}{2}, \quad \text { and } |r_{k}-r_{k-1}|=(R_0-R)2^{-k}
\]
By (\ref{23}), we have
\begin{eqnarray}\label{24}
||S||_{L^{q_{k+1}+\frac{1}{2}}(B_{r_{k+1}}(z))}&\leq& \Big[Cq_k(1+\frac{1}{(r_{k+1}-r_k)^2}\Big]^{\frac{1}{q_k}}||S||^{a_k}_{L^{q_k+\frac{1}{2}}(B_{r_k}(z))}\\\nonumber
&\leq& q_k^{\frac{1}{q_k}}\Big(C(1+\frac{1}{(R_0-R)^2})\Big)^{\frac{1}{q_k}}2^{\frac{2k}{q_k}}||S||^{a_k}_{L^{q_k+\frac{1}{2}}(B_{r_k}(z))}.
\end{eqnarray}
where $a_k:=\frac{q_k+\frac{1}{2}}{q_k}$. By iteration, it follows from (\ref{24}) that
\begin{eqnarray}\label{25}
&&||S||_{L^{q_{k+1}+\frac{1}{2}}(B_{r_{k+1}}(z))}\\\nonumber
&\leq&\Big[ \prod^k_{i=1}q_i^{\frac{1}{q_i}}\Big(C(1+\frac{1}{(R_0-R)^2})\Big)^{\frac{1}{q_i}}2^{\frac{2i}{q_i}}\Big]^{\prod^k_{i=1}a_i}||S||^{\prod^k_{i=1}a_i}_{L^{q_1+\frac{1}{2}}(B_{r_1}(z))}.
\end{eqnarray}
Notice that $a_k=\frac{q_k+\frac{1}{2}}{q_k}=\frac{\frac{q_{k-1}m}{m-1}}{q_k}=\frac{m}{m-1}\frac{q_{k-1}}{q_k}$, so
\[
\prod^k_{i=1}a_i=\Big(\frac{m}{m-1}\Big)^k \frac{q_0}{q_1}\cdots \frac{q_{k-1}}{q_k}=\Big(\frac{m}{m-1}\Big)^k  \frac{q_0}{q_k}
\]
and thus
\[
\lim_{k\rightarrow \infty}\prod^k_{i=1}a_i=q_0=\frac{m+1}{2}.
\]
Moreover,
 \[
 \prod^k_{i=1}q_i^{\frac{1}{q_i}}\Big(C(1+\frac{1}{(R_0-R)^2})\Big)^{\frac{1}{q_i}}2^{\frac{2i}{q_i}}=\prod^k_{i=1}q_i^{\frac{1}{q_i}}
 \Big(C(1+\frac{1}{(R_0-R)^2})\Big)^{\sum_{i=1}^k\frac{1}{q_i}} 2^{\sum_{i=1}^k\frac{2i}{q_i}}.
\]
When $k\rightarrow\infty$, it is easy to show that $\sum_{i=1}^{\infty}\frac{1}{q_i}<\infty \text{ and } \sum_{i=1}^{\infty}\frac{2i}{q_i}<\infty$. Notice also that $\log(\prod^{\infty}_{i=1}q_i^{\frac{1}{q_i}})<\infty$. Thus,
\[
\lim_{k\rightarrow\infty} \prod^k_{i=1}q_i^{\frac{1}{q_i}}\Big(C(1+\frac{1}{(R_0-R)^2})\Big)^{\frac{1}{q_i}}2^{\frac{2i}{q_i}}<\infty.
\]
It follows from (\ref{25}), by letting $k\rightarrow \infty$,
\begin{equation}\label{26}
||S||_{L^{\infty}}\leq C ||S||^{\frac{m+1}{2}}_{L^{q_1+\frac{1}{2}}(B_{R_0}(z))}.
\end{equation}
Choosing now $\s=q_1+\frac{1}{2}=\frac{m}{m-1}+\frac{m}{2}$ in (\ref{17}), we finally obtain
\begin{equation}\label{27}
||S||_{L^{\infty}}\leq C,
 \end{equation}
where $C$ is a positive constant depending on $ K,  |d\o|_{\o},  |R|_{\o}, |\n R|_{\o}, |T|_{\o}, |\n T|_{\o},  dist (\O', \partial \O )$ and $|\n ^s f|_{\o}$ , $s=0, 1, 2, 3$.

\vspace{1cm}
\medskip
\section{Appendix}
\par
As mentioned in the introduction, using the idea from \cite{PSS}, we give a new proof for the elliptic inequality (\ref{elliptic}) in this section.
\medskip
\par
\noindent{\bf Proof of the elliptic inequality (\ref{elliptic}):}
\par
Let $\n$ and $\tn$ denote the Chern connections corresponding to the Hermitian metircs $\o$ and $\o+\sqrt{-1}\p\bar\p\phi$ respectively. Define
\begin{equation}\label{00}
h=\tg\cdot g^{-1}
\end{equation}
and
\[
h^j_i=\tg_{i\bar k}g^{j\bar k}, \quad (h^{-1})^j_i =g_{i\bar k}\tg^{j\bar k}.
\]
In fact, $h$ can be thought to be an endomorphism $h: T^{1, 0}(M)\longrightarrow T^{1, 0}(M)$, such that $\tg(X, Y)= g(h(X), Y)$.
\par
Set
\begin{equation}\label{01}
S=\tg^{j\bar r}\tg^{s\bar k}\tg^{m\bar l}\phi_{j\bar k m}\phi_{\bar r s\bar l},
\end{equation}
where $\phi_{j\bar k m}=\n_m\n_{\bar k}\n_j \phi$.
\par
By (\ref{00}), we have
\begin{eqnarray}\label{02}
\td\th &=& \p \tg\cdot \tg^{-1}=\p (h \cdot g) \cdot g^{-1} h^{-1} \\\nonumber
&=& \p h \cdot g\cdot g^{-1}\cdot h^{-1}+h\cdot \p g\cdot g^{-1}\cdot h^{-1}\\\nonumber
&=&\p h\cdot h^{-1}+h\cdot \th\cdot h^{-1}\\\nonumber
&=& \p h\cdot h^{-1}+h\cdot \th \cdot h^{-1}- \th \cdot h\cdot h^{-1}+\th\\\nonumber
&=&\th +(\n^{1, 0}h)\cdot h^{-1}.
\end{eqnarray}
\begin{eqnarray}\label{03}
\td R&=&\bar \p \td\th=\bar\p (\th +(\n^{1, 0}h)\cdot h^{-1})\\\nonumber
&=& R+\bar\p ((\n^{1, 0}h)\cdot h^{-1}).
\end{eqnarray}
By similar computation, we can get
\begin{equation}\label{04}
\th=\p g\cdot g^{-1}= \td\th- h^{-1}(\tn^{1, 0}h),
\end{equation}
\begin{equation}\label{05}
R=\td R-\bar\p(h^{-1}\cdot (\tn^{1, 0}h)).
\end{equation}
Now, using the definitions, one can see that
\begin{equation*}
\phi_{j\bar k m}=(\n_m \tg)(\p_j, \bar \p_k)=\tg_{j\bar k, m}.
\end{equation*}
Thus,
\begin{equation}\label{06}
S=\tg^{j\bar r}\tg^{s\bar k}\tg^{m\bar l}\phi_{j\bar k m}\phi_{\bar r s\bar l}=|\n^{1, 0}\tg|_{\tg}^2.
\end{equation}
On the other hand,
\begin{eqnarray*}
\n_m \tg=\n_m(h\cdot g)=\n_m h\cdot g=\Big(\frac{\p}{\p z^m}h+h\cdot \th_m-\th_m\cdot h\Big)\cdot g,
\end{eqnarray*}
so
\begin{eqnarray*}
\td\n_m h&=&\frac{\p}{\p z^m}h+h\cdot \td\th_m-\td\th_m\cdot h\\
&=&\frac{\p}{\p z^m}h+h\cdot \th_m-\th_m\cdot h+h\cdot(\n_m h)\cdot h^{-1}-\n_m h\\
&=& h\cdot (\n_m h)\cdot h^{-1}.
\end{eqnarray*}
Thus,
\begin{equation*}
\n_m \tg=\n_m h\cdot g=h^{-1}\cdot (\td\n_m h)\cdot h\cdot g=h^{-1}\cdot (\td\n_m h)\cdot \tg.
\end{equation*}
Finally we end up with the formula
\begin{equation}\label{07}
S=|\n^{1, 0}\tg|_{\tg}^2=|h^{-1}\cdot (\td\n^{1, 0}h)|_{\tg}^2=|\td\th-\th|_{\tg}^2
\end{equation}
i.e. $S$ can be thought as the $\tg$-norm of the difference between the two connection 1-forms.
\par
Now, we can deduce the elliptic inequality:
\begin{eqnarray}\label{08}
\td\tri S&=& \td\tri |h^{-1}\cdot (\td \n^{1, 0}h)|^2_{\tg}\\\nonumber
&=&\tg^{i\bj}\p_i\p_{\bj}<h^{-1}\cdot(\td\n^{1, 0}h) ,  \overline{h^{-1}\cdot(\td\n^{1, 0}h)}>_{\tg}\\\nonumber
&=&\tg^{i\bj}\p_i\Big(<\td\n_{\bj}(h^{-1}\cdot(\td\n^{1, 0}h) ),  \overline{h^{-1}\cdot(\td\n^{1, 0}h)}>_{\tg}\\\nonumber
&&+<(h^{-1}\cdot(\td\n^{1, 0}h) ),  \overline{\td\n_{j}h^{-1}\cdot(\td\n^{1, 0}h)}>_{\tg}\Big)\\\nonumber
&=&\tg^{i\bj}<\td\n_{i}\td\n_{\bj}(h^{-1}\cdot(\td\n^{1, 0}h) ), \overline{h^{-1}\cdot (\td\n^{1, 0}h)}>_{\tg}\\\nonumber
&&+\tg^{i\bj}<h^{-1}\cdot (\td\n^{1, 0}h),  \overline{\td\n_{\bi}\td\n_{j}(h^{-1}\cdot(\td\n^{1, 0}h)})>_{\tg}\\\nonumber
&&+|\td\n^{1, 0}(h^{-1}\cdot (\td\n^{1, 0}h))|_{\tg}^2+|\td\n^{0, 1}(h^{-1}\cdot (\td\n^{1, 0}h))|_{\tg}^2.
\end{eqnarray}
Using the relation $R=\td R-\bar \p(h^{-1}\cdot(\td\n^{1, 0}h))$, we have
\begin{equation}\label{09}
\tg^{i\bj}\td\n_i\td\n_{\bj}(h^{-1}\cdot (\td\n_t^{1, 0}h))^l_m=\tg^{i\bj}\td\n_i\Big(\td R^l_{m t\bj}-R^l_{m t\bj}\Big).
\end{equation}

\medskip
\par
Recall the Bianchi identities of curvature forms which can be found in \cite{KN}(p. 135):
\begin{eqnarray}
\label{010}&&\sum ( R(X, Y) Z)= \sum T(T(X, Y), Z)+(\n_X T)(Y, Z);\\
\label{011}&& \sum\{\n_X R(Y, Z)+R(T(X, Y), Z)\}=0,
\end{eqnarray}
where $X, Y, Z\in TM$ and T is the torsion of the connection $\n$ (recall that $\n$ is not necessarily the Levi-Civita connection), while $\sum$ denotes the cyclic sum with respect to $X, Y, Z$.

\medskip
\par
By the first Bianchi identity (\ref{010}), one obtains
\begin{eqnarray*}
&&\td R(\p_i, \p_{\bj})\p_{m}+\td R(\p_{\bj}, \p_{m})\p_{i}+\td R(\p_{m}, \p_{i})\p_{\bj}\\
&=&\td T\Big(\td T(\p_{i}, \p_{\bj}), \p_{m}\Big)+\td T\Big(\td T(\p_{\bj}, \p_{m}), \p_{i}\Big)+\td T\Big(\td T(\p_{m}, \p_{i}), \p_{\bj}\Big)\\
&&+(\td\n_{i}\td T)(\p_{\bj}, \p_{m})+(\td\n_{\bj}\td T)(\p_{m}, \p_{i})+(\td\n_{m}\td T)(\p_{i}, \p_{\bj}).
\end{eqnarray*}
Recall the fact that $\td R^{2, 0}=\td R^{0, 2}=0$, $\td T^{1, 1}=0$ (since $\td\n$ is the Chern connection) and $\td T(\p_m, \p_i)\in T^{1, 0}(M)$. Also
\begin{eqnarray*}
&&\td T(\p_{i}, \p_{\bj})=\td T(\p_{\bj}, \p_{m})=(\td\n_{i}\td T)(\p_{\bj}, \p_{m})=(\td\n_{m}\td T)(\p_{i}, \p_{\bj})=0,\\
&& \td R(\p_{m}, \p_{i})\p_{\bj}=0.
\end{eqnarray*}
Thus,
\begin{equation*}
\td R(\p_i, \p_{\bj})\p_{m}+\td R(\p_{\bj}, \p_{m})\p_{i}=(\td\n_{\bj}\td T)(\p_{m}, \p_{i}).
\end{equation*}
By definition $\td R(\p_i, \p_{\bj})\p_{m}=\td R^l_{m i\bj}\p_l$ and $\td R^l_{m i\bj}=-\td R^l_{m \bj i}$, so we get
\begin{equation}\label{012}
\td R^l_{m i\bj}=\td R^l_{i m \bj}+\td T^l_{m i, \bj}.
\end{equation}
Similarly, one can also obtain
\begin{equation}\label{013}
\td R^{\bar l}_{\bar k i\bj}=\td R^{\bar l}_{\bj i \bar k}+\td T^{\bar l}_{\bj \bar k, i}.
\end{equation}

Moreover, by the second Bianchi identity (\ref{011}) and following the same step as above we have
\begin{eqnarray*}
\td R^l_{m t\bj, i}+\td R^l_{m\bj i, t}+\td R^l_{m i t, \bj}=-\td R(\td T(\p_i, \p_t), \p_{\bj})-\td R(\td T(\p_t, \p_{\bj}), \p_{i})-\td R(\td T(\p_{\bj}, \p_i), \p_{t})
\end{eqnarray*}
and $\td R^l_{m i t, \bj}=0,  \td T(\p_t, \p_{\bj})=\td T(\p_{\bj}, \p_i)=0$.
Thus,
\begin{equation}\label{014}
\td R^l_{m i\bj, t}=\td R^l_{m t\bj, i}+\td T^s_{i t}\td R^l_{m s\bj}.
\end{equation}
Now, using the identities (\ref{012}), (\ref{013}) and (\ref{014}), we obtain
\begin{eqnarray}\label{015}
\tg^{i\bj}\td\n_i \td R^l_{m t\bj}&=&\tg^{i\bj}\td R^l_{m t\bj, i}=\tg^{i\bj}\td R^l_{m i\bj, t}-\tg^{i\bj}\td T^s_{i t}\td R^{l}_{m s\bj}\\\nonumber
&=&\tg^{i\bj}\td R_{m \bar k i\bj, t}\tg^{l\bar k}-\tg^{i\bj}\td T^s_{i t}\td R^l_{m s\bj}\\\nonumber
&=&\tg^{i\bj}(\td R_{i\bar k m\bj, t}+\td T^s_{m i, \bj t}\tg_{s\bar k})\tg^{l\bar k}-\tg^{i\bj}\td T^s_{i t}\td R^l_{m s\bj}\\\nonumber
&=&-\tg^{i\bj}\td R_{\bar k i m \bj, t}\tg^{l\bar k}+\tg^{i\bj}\td T^l_{m i, \bj t}-\tg^{i\bj}\td T^s_{i t}\td R^l_{m s\bj}\\\nonumber
&=&-\tg^{i\bj}\td R_{\bj i m \bar k, t}\tg^{l\bar k}-\tg^{i\bj}\td T^{\bar l}_{\bj \bar k, m t}\tg_{i\bar l}\tg^{l\bar k}+\tg^{i\bj}\td T^l_{m i, \bj t}-\tg^{i\bj}\td T^s_{i t}\td R^l_{m s\bj}\\\nonumber
&=&\tg^{i\bj}\td R_{i\bj m \bar k, t}\tg^{l\bar k}-\tg^{i\bj}\td T^{\bar l}_{\bj \bar k, m t}\tg_{i\bar l}\tg^{l\bar k}+\tg^{i\bj}\td T^l_{m i, \bj t}-\tg^{i\bj}\td T^s_{i t}\td R^l_{m s \bj}\\\nonumber
&=& \td R^i_{i m \bar k, t}\tg^{l\bar k}-\tg^{i\bj}\td T^{\bar l}_{\bj \bar k, m t}\tg_{i\bar l}\tg^{l\bar k}+\tg^{i\bj}\td T^l_{m i, \bj t}-\tg^{i\bj}\td T^s_{i t}\td R^l_{m s \bj}
\end{eqnarray}
From the Monge-Amp\`ere equation (\ref{eqn}), it follows that
\begin{equation}\label{016}
\td R^i_{i m \bar k, t}=\td \n_t R^i_{i m \bar k}-\td \n_t f_{m \bar k}.
\end{equation}
In the following, we denote $\epsilon =O(S^{\alpha })$ if there is a constant C depending only on $ K,  |d\o|_{\o},  |R|_{\o}, |\n R|_{\o}, |T|_{\o}, |\n T|_{\o}$ and $|\n ^s f|_{\o}$ , $s=0, 1, 2, 3$, such that $\epsilon \leq C S^{\alpha }$. Note that $\td\n$ is $O(S^{\frac{1}{2}})$, so
\begin{equation}\label{017}
\td R^i_{i m \bar k, t}\tg^{l\bar k}=O(S^{\frac{1}{2}})+O(1).
\end{equation}
For the second term in (\ref{015})
\begin{eqnarray}\label{018}
\td T^{\bar s}_{\bj \bar k, m t}&=&\Big((\p_{\bj}g_{n \bar k}-\p_{\bar k}g_{n\bj})\tg^{n\bar s}\Big)_{m t}\\\nonumber
&=&(T_{\bj \bar k n} \tg^{n\bar s})_{m t}= \td\n_t\td \n_m T_{\bj \bar k n}\tg^{n\bar s}\\\nonumber
&=&\td \n_t(\n_m T_{\bj\bar k n}-(\td\th_m-\th_m)^l_n T_{\bj \bar k l})\tg^{n\bar s}\\\nonumber
&=&\Big(\n_t (\n_m T_{\bj \bar k n})-(\td\th_t-\th_t)^l_m\n_l T_{\bj\bar k n}-\td\n_t((\td\th_m-\th_m)^l_n)T_{\bj\bar k l}\\\nonumber
&&-(\td\th_t-\th_t)^l_n \n_m T_{\bj\bar k l}-(\td\th_m-\th_m)^l_n(\n_t T_{\bj \bar k l}-(\td\th_t-\th_t)^s_l T_{\bj \bar k s})\Big)\tg^{n\bar s}.\\\nonumber
\end{eqnarray}
Again, by the fact that $\td\n$ is $O(S^{\frac{1}{2}})$ and $|h^{-1}\cdot (\td\n^{1, 0}h)|_{\tg}$ is also $O(S^{\frac{1}{2}})$, we have
\begin{equation}\label{019}
|\tg^{i\bj}\td T^{\bar l}_{\bj \bar k, m t}\tg_{i\bar l}\tg^{l\bar k}|\leq O(S^{\frac{1}{2}})+O(S)+C |\td \n^{1, 0}(h^{-1}\cdot (\td\n^{1, 0}h))|+O(1).
\end{equation}

Similarly, we can get the estimate for the last two terms in (\ref{015})
\begin{eqnarray}
\label{019}|\tg^{i\bj}\td T^l_{m i, \bj t}|&\leq& O(S^{\frac{1}{2}})+O(S)+C |\td \n^{0, 1}(h^{-1}\cdot (\td\n^{1, 0}h))|+O(1),\\
\label{020}|\tg^{i\bj}\td T^s_{i t}\td R^l_{m s \bj}|&\leq&C |\td \n^{0, 1}(h^{-1}\cdot (\td\n^{1, 0}h))|+O(1).
\end{eqnarray}
Put the above estimates (\ref{015})-(\ref{020}) into (\ref{09}), we can conclude that
\begin{eqnarray}\label{021}
&&|\tg^{i\bj}\td\n_i\td\n_{\bj}(h^{-1}\cdot (\td\n_t^{1, 0}h))^l_m|\\\nonumber&\leq& O(S^{\frac{1}{2}})+O(S)
+C |\td \n^{1, 0}(h^{-1}\cdot (\td\n^{1, 0}h))|+C |\td \n^{0, 1}(h^{-1}\cdot (\td\n^{1, 0}h))|.
\end{eqnarray}
One the other hand,
\begin{equation*}
\tg^{i\bj}\td\n_{\bi}\td\n_{j}(h^{-1}\cdot (\td\n^{1, 0}h))=\tg^{i\bj}\td\n_{j}\td\n_{\bi}(h^{-1}\cdot (\td\n^{1, 0}h))-(\tg^{i\bj}\td R^l_{m i\bj})\#(h^{-1}\cdot (\td\n^{1, 0}h))
\end{equation*}
where
\begin{equation*}
\begin{array}{lll}
&&(\tg^{i\bj}\td R^l_{m i\bj})\#(h^{-1}\cdot (\td\n^{1, 0}h))\\
&=&\tg^{i\bj}\{h^{-1}\cdot (\td\n^{1, 0}_{t}h)_{m}^{s}\td R^l_{s i\bj}-h^{-1}\cdot (\td\n^{1, 0}_{s}h)_{m}^{l}\td R^s_{t i\bj}-h^{-1}\cdot (\td\n^{1, 0}_{t}h)_{s}^{l}\td R^s_{m i\bj}\}\\
&& dz^{t}\otimes dz^{m}\otimes \frac{\partial }{\partial z^{l}}\\
\end{array}
\end{equation*}
and
\begin{eqnarray*}
\tg^{i\bj}\td R^l_{m i\bj}&=&\tg^{i\bj}\td R^l_{i m\bj}+\tg^{i\bj}\td T^l_{m i, \bj}=\tg^{i\bj}\td R_{i\bj m \bar k}\tg^{l\bar k}+\tg^{i\bj}\td T^{\bar s}_{\bj \bar k, m}\tg_{i\bar s}\tg^{l\bar k}+\tg^{i\bj}\td T^l_{m i, \bj}.
\end{eqnarray*}
Thus
\begin{eqnarray*}
 |\tg^{i\bj}\td R^l_{m i\bj}|&\leq& O(S^{\frac{1}{2}}) +O(1).
\end{eqnarray*}
Hence we conclude that
\begin{eqnarray}\label{022}
&&|\tg^{i\bj}\td\n_{\bi}\td\n_{j}(h^{-1}\cdot (\td\n^{1, 0}h))|\\\nonumber
&\leq& |\tg^{i\bj}\td\n_{j}\td\n_{\bi}(h^{-1}\cdot (\td\n^{1, 0}h))|+|(\tg^{i\bj}\td R^l_{m i\bj})\#(h^{-1}\cdot (\td\n^{1, 0}h))|\\\nonumber
&\leq&O(S^{\frac{1}{2}})+O(S)+C |\td \n^{1, 0}(h^{-1}\cdot (\td\n^{1, 0}h))|+C |\td \n^{0, 1}(h^{-1}\cdot (\td\n^{1, 0}h))|.
\end{eqnarray}
\par
Finally, by (\ref{08}) and (\ref{021}), (\ref{022}), we obtain the elliptic inequality:
\begin{equation}\label{023}
\td\tri S\geq -C_1 S^{\frac{3}{2}}-C_2
\end{equation}
where $C_1, C_2$ are positive constants depending only on $ K,  |d\o|_{\o},  |R|_{\o}, |\n R|_{\o}, |T|_{\o}, |\n T|_{\o}$ and $|\n ^s f|_{\o}$ , $s=0, 1, 2, 3$.
\medskip

{\bf Acknowledgements:} The authors would like to thank Prof. Pengfei Guan and Slawomir Dinew for the numerous helpful discussions on this problem. The note was written while the first named author was visiting McGill University. He would like to thank this institution for the hospitality.

\end{document}